\documentclass[10pt,journal]{IEEEtran}
\usepackage{times}
\usepackage{amsmath}
\usepackage{amsfonts}
\usepackage[dvips]{graphicx}
\newcommand{\Vector}[1]{{\bf \lowercase{#1}}}
\newcommand{\Matrix}[1]{{\bf \uppercase{#1}}}

\renewcommand{\baselinestretch}{1.15}

\title{Generalized Canonical Correlation Analysis and Its Application to Blind Source Separation Based on a Dual-Linear Predictor Structure}

\author{Wei Liu\\ \vspace*{0.4cm}
       {\normalsize Communications Research Group}\\ 
       {\normalsize Dept. of Electronic \& Electrical Engineering}\\
        {\normalsize University of Sheffield, UK}\\
         { \tt w.liu@sheffield.ac.uk }}

\begin{document}
\maketitle
\begin{abstract}
Blind source separation (BSS) is one of the most important and established research topics
in signal processing and many algorithms have been proposed based on different statistical
properties of the source signals. For second-order statistics (SOS) based methods, canonical
correlation analysis (CCA) has been proved to be an effective solution to the problem. In this
work, the CCA approach is generalized to accommodate the case with added white noise
and it is then applied to the BSS problem for noisy mixtures. In this approach, the noise
component is assumed to be spatially and temporally white, but the variance information of
noise is not required. An adaptive blind source extraction algorithm is derived based on this
idea and a further extension is proposed by employing a dual-linear predictor structure for
blind source extraction (BSE).
\end{abstract}




%
%
\section{INTRODUCTION}
The problem of blind source separation (BSS) has been studied extensively in the past and a plethora of algorithms have been proposed in the past based on statistical properties of the source signals~\cite{hyvarinen01a,cichocki03a}. For those based on the second-order statistics (SOS), one particular class of them is those based on the canonical correlation analysis (CCA) approach~\cite{schell95a,galy00a,borga01a,friman02a,liu06d,liu07l,li09a,liu13b,liu13l}, where the demixing matrix is found by maximizing the autocorrelation of each of the recovered signals. This approach rests on the idea that the sum of any uncorrelated signals has an autocorrelation whose value is less or equal to the maximum value of individual signals. As shown in \cite{liu07l}, the maximization of the autocorrelation value is equivalent to finding the generalised eigenvectors within the matrix pencil
approach~\cite{chang98a}.

As shown in ~\cite{liu07l}, the CCA approach will work for the noise-free situations. For noisy mixtures, its performance will no doubt degrade. If we can estimate the variance of white noise in the mixtures, then we can remove the effect of noise from the mixtures before applying the noise-free BSS algorithms~\cite{liu06e,liu06f}. However, the variance of noise in the mixtures is not always available and in this case the algorithm proposed in ~\cite{liu06f} will not work.

In this paper, we will generalise the traditional CCA to include the case with noise and then apply it to the separation problem of noisy mixtures~\cite{liu08h,liu08g}. A key advantage of this approach is that successful separation of source signals can be achieved  without estimation of the white noise parameters. An online adaptive algorithm is derived accordingly as its adaptive realisation.

Moreover, similar to the CCA case \cite{liu07l}, we have also related the GCCA approach to a dual-linear predictor based blind source extraction (BSE) structure, and an adaptive algorithm based on such a structure is also derived with rigorous proof for its effectiveness in this context.

This paper is organised as follows. In
Section~\ref{sec:GCCA}, the generalised CCA approach will be provided with a detailed proof and analysis about its condition on which it can be applied to the BSS problem. The class of adaptive BSE algorithms  is derived
in Sec.~\ref{sec:GCCA_adaptive}. Simulation results are shown in
Section~\ref{sec:GCCA_sim} and conclusions drawn in
Section~\ref{sec:GCCA_concl}.
%
%
\section{Generalised CCA and Its Application to BSS}
\label{sec:GCCA}

\subsection{Overview of the CCA Approach to BSS}
\label{sec:GCCA_cca}

The instantaneous mixing problem in BSS with $M$ mixtures, $L$ sources and a mixing matrix $\Matrix{A}$ can be expressed as
\begin{equation}
     \Vector{x}[n] = \Matrix{A}\cdot\Vector{s}[n],
\label{eqn:intro_instan}
\end{equation}
with
\begin{eqnarray}
     \Vector{s}[n] &=& \left[s_0[n] \; s_1[n] \; \cdots \; s_{L\!-\!1}[n]\right]^\mathrm{T} \nonumber \\
      \Vector{x}[n]&=& \left[ x_0[n] \; x_1[n]\; \cdots \; x_{M\!-\!1}[n]\right]^\mathrm{T} \nonumber \\
     \left[\Matrix{A}\right]_{m,l} & = & a_{m,l} \;, m=0,\dots,
     M-1,\; l=0,\dots, L-1.
\end{eqnarray}
For BSS employing second-order statistics (SOS), we assume the sources are spatially uncorrelated and have different temporal
structures:
\begin{eqnarray}
     \Matrix{R}_{ss}[0]&=&E\{\Vector{s}[n]\Vector{s}^T[n]\} = \mbox{diag}\{\rho_0[0],\rho_1[0], \dots,
     \rho_{L-1}[0]\}\nonumber \\
     \Matrix{R}_{ss}[\Delta n]&=&E\{\Vector{s}[n]\Vector{s}^T[n-\Delta n]\}= \mbox{diag}\{\rho_0[\Delta n],\rho_1[\Delta n], \dots, \rho_{L-1}[\Delta n]\}
     \label{eqn:gcca_sourcecorrelation}
\end{eqnarray}
with $\rho_l[n]$ being the autocorrelation function of the $l$th source signal and $\rho_l[\Delta n]\neq 0$ for some nonzero delays $\Delta n$.

The BSS problem can be solved in one single step by the CCA approach. In CCA~\cite{anderson84a}, two sets of zero-mean variables $\Vector{z}_1[n]$ (with $q_1$ components) and $\Vector{z}_2[n]$ (with $q_2$ components) with a joint distribution are
considered. For convenience, we assume $q_1\leq q_2$. The linear combination of the variables in each
of the vectors is respectively given by
\begin{eqnarray}
     a_0[n] &=& \boldsymbol{\alpha}_0^T\Vector{z}_1[n] \nonumber \\
     b_0[n] &=& \boldsymbol{\beta}_0^T\Vector{z}_2 [n]\;,
\end{eqnarray}
where $\boldsymbol{\alpha}_0$ and $\boldsymbol{\beta}_0$ are
vectors containing the combination coefficients and they are determined by maximizing the correlation between $a_0$
and $b_0$
\begin{equation}
    \max_{\boldsymbol{\alpha}_0,\boldsymbol{\beta}_0} J_0(\boldsymbol{\alpha}_0,\boldsymbol{\beta}_0)\;
    \label{eqn:cca_firstcorrelation}
\end{equation}
with
\begin{eqnarray}
J_0(\boldsymbol{\alpha}_0,\boldsymbol{\beta}_0)&=&\frac{E\{a_0[n]b_0[n]\}}{\sqrt{E\{a_0^2[n]\}E\{b_0^2[n]\}}}\nonumber \\
&=&\frac{\boldsymbol{\alpha}_0^T\boldsymbol{\Sigma}_{12}\boldsymbol{\beta}_0}{\sqrt{(\boldsymbol{\alpha}_0^T\boldsymbol{\Sigma}_{11}[0]\boldsymbol{\alpha}_0)(\boldsymbol{\beta}_0^T\boldsymbol{\Sigma}_{22}[0]\boldsymbol{\beta}_0)}}\;,
\label{eqn:cca_firstcorrelation1}
\end{eqnarray}
where $\boldsymbol{\Sigma}_{11}[0]=E\{\Vector{z}_1[n]\Vector{z}_1^T[n]\}$,
$\boldsymbol{\Sigma}_{12}=E\{\Vector{z}_1[n]\Vector{z}_2^T[n]\}$,
$\boldsymbol{\Sigma}_{22}=E\{\Vector{z}_2[n]\Vector{z}_2^T[n]\}$ and $E\{\cdot\}$ denotes the statistical expectation operator.

After finding the first pair of optimal vectors
$\boldsymbol{\alpha}_0$ and $\boldsymbol{\beta}_0$, we can proceed
to find the second pair $\boldsymbol{\alpha}_1$ and
$\boldsymbol{\beta}_1$ which maximizes the correlation and at the
same time ensures that the new pair of combinations
$\{a_1[n],\;b_1[n]\}$ is uncorrelated with the first set
$\{a_0[n],\;b_0[n]\}$. This process is repeated until we find all the
$\min(q_1,q_2)=q_1$ pairs of optimal vectors
$\boldsymbol{\alpha}_i$ and $\boldsymbol{\beta}_i$, $i=0$, $\dots$,
$q_1-1$.

 It has been shown that  $\boldsymbol{\alpha}_i$  can be obtained by solving the following generalized eigenvalue
problem~\cite{liu07l}
\begin{equation}
 \boldsymbol{\Sigma}_{12}\boldsymbol{\Sigma}_{22}^{-1}\boldsymbol{\Sigma}_{21}\boldsymbol{\alpha}_i=\lambda_i^2\boldsymbol{\Sigma}_{11}\boldsymbol{\alpha}_i\;.
 \label{eqn:cca_eigenvector}
\end{equation}
 $\boldsymbol{\beta}_i$ can be found in the same way by changing the subscripts
of the matrices in \eqref{eqn:cca_eigenvector} accordingly.

 To apply CCA to the BSS problem~\cite{liu07l}, we choose the vector $\Vector{x}[n]$ as $\Vector{z}_1$ in
CCA and $\Vector{x}[n-\Delta_n]$ as $\Vector{z}_2$. Then the eigenvalue problem in \eqref{eqn:cca_eigenvector} becomes
\begin{equation}
 \Matrix{R}_{xx}[\Delta n]\Matrix{R}_{xx}[0]^{-1}\Matrix{R}_{xx}[\Delta n]\boldsymbol{\alpha}_i=\lambda_i^2\Matrix{R}_{xx}[0]\boldsymbol{\alpha}_i\;.
 \label{eqn:cca_bsseigenvector1}
\end{equation}

In the context of BSS, $\boldsymbol{\alpha}_i$ and $\boldsymbol{\beta}_i$ are
the same and we use  $\Vector{w}_i$ to represent it.
\eqref{eqn:cca_bsseigenvector1} can be simplified as
\begin{equation}
 \Matrix{R}_{xx}[0]^{-1}\Matrix{R}_{xx}[\Delta
 n] \Vector{w}_i=\lambda_i \Vector{w}_i\;.
 \label{eqn:cca_bsseigenvector3}
\end{equation}
Multiplying both sides with $\Matrix{R}_{xx}[0]$, we arrive at the following generalised eigenvector problem
\begin{equation}
 \Matrix{R}_{xx}[\Delta
 n]\Vector{w}_i=\lambda_i\Matrix{R}_{xx}[0]\Vector{w}_i\;.
 \label{eqn:cca_bsseigenvector4}
\end{equation}
Moreover, the correlation maximization problem in \eqref{eqn:cca_firstcorrelation} becomes
\begin{equation}
    \max_{\Vector{w}_0} J_0(\Vector{w}_0)=\frac{\Vector{w}_0^T\Matrix{R}_{xx}[\Delta
n]\Vector{w}_0}{\Vector{w}_0^T\Matrix{R}_{xx}[0]\Vector{w}_0}=\frac{\Vector{w}_0^T\Matrix{A}\Matrix{R}_{ss}[\Delta
n]\Matrix{A}^T\Vector{w}_0}{\Vector{w}_0^T\Matrix{A}\Matrix{R}_{ss}[0]\Matrix{A}^T\Vector{w}_0}\;,
    \label{eqn:cca_bssfirstcorrelation}
\end{equation}
and we can prove that by CCA the source signals will be recovered completely~\cite{liu07l}.
However, with added noise, the proof given in the noise-free case will not be valid any more since the denominator in \eqref{eqn:cca_firstcorrelation1} will have a noise component. As a result, the performance of the CCA approach will degrade with increasing noise level.

\subsection{Generalised CCA (GCCA)}
\label{sec:GCCA_gcca}

For noisy mixtures, $\Vector{x}[n]$ is given by
\begin{equation}
     \Vector{x}[n] = \Matrix{A}\Vector{s}[n]+\Vector{v}[n],
\label{eqn:gcca_noisymixtures}
\end{equation}
where $\Vector{v}[n]$ is the additive noise vector, which is spatially and temporally white and uncorrelated with
the source signals. Its correlation matrix is given by
\begin{equation}
\Matrix{R}_{vv}[\Delta_n]= E\{\Vector{v}[n]\Vector{v}^T[n-\Delta_n]\}=\left\{ \begin{array}{ll}
        \Matrix{0} & \mbox{for }   \Delta_n\neq 0 \\
        \sigma_v^2\Matrix{I} &\mbox{for } \Delta_n=0
    \end{array} \right. .
    \label{eqn:gcca_noisecorrelation}
\end{equation}
where $\Matrix{I}$ is the identity
matrix and $\sigma_v^2$ is the variance of noise.

Similarly, we can form a modified CCA problem for two set of variables with added white noise. Now consider the two sets of zero-mean variables
\begin{eqnarray}
\hat{\Vector{z}}_1[n]&=&\Vector{z}_1[n]+\Vector{v}_1[n]\nonumber \\
\hat{\Vector{z}}_2[n]&=&\Vector{z}_2[n]+\Vector{v}_2[n]
\end{eqnarray}
and their corresponding linear combinations:
\begin{eqnarray}
     \hat{a}_0[n] &=& \boldsymbol{\alpha}_0^T\hat{\Vector{z}}_1[n] \nonumber \\
     \hat{b}_0[n] &=& \boldsymbol{\beta}_0^T\hat{\Vector{z}}_2 [n]\;,
\end{eqnarray}
where $\Vector{v}_1$ and $\Vector{v}_2$ are the added white noise vectors and not correlated with each other. Now the two vectors $\boldsymbol{\alpha}_0$ and $\boldsymbol{\beta}_0$ are not given by \eqref{eqn:cca_firstcorrelation} any more, but by
\begin{equation}
    \max_{\boldsymbol{\alpha}_0,\boldsymbol{\beta}_0} \hat{J}_0(\boldsymbol{\alpha}_0,\boldsymbol{\beta}_0)\;
    \label{eqn:gcca_firstcorrelation}
\end{equation}
with
\begin{eqnarray}
\hat{J}_0(\boldsymbol{\alpha}_0,\boldsymbol{\beta}_0)&=&\frac{E\{\hat{a}_0[n]\hat{b}_0[n]\}}{E\{\hat{a}_0[n]\hat{a}_0[n-\Delta_0]\}E\{\hat{b}_0[n]\hat{b}_0[n-\Delta_0]\}}\nonumber \\
&=&\frac{\boldsymbol{\alpha}_0^T\hat{\boldsymbol{\Sigma}}_{12}\boldsymbol{\beta}_0}{\sqrt{(\boldsymbol{\alpha}_0^T\hat{\boldsymbol{\Sigma}}_{11}[\Delta_0]\boldsymbol{\alpha}_0)(\boldsymbol{\beta}_0^T\hat{\boldsymbol{\Sigma}}_{22}[\Delta_0]\boldsymbol{\beta}_0)}}\;,
\label{eqn:gcca_firstcorrelation1}
\end{eqnarray}
where $\hat{\boldsymbol{\Sigma}}_{12}=E\{\hat{\Vector{z}}_1[n]\hat{\Vector{z}}_2^T[n]\}$, $\hat{\boldsymbol{\Sigma}}_{11}[\Delta_0]=E\{\hat{\Vector{z}}_1[n]\hat{\Vector{z}}_1^T[n-\Delta_0]\}$,
$\hat{\boldsymbol{\Sigma}}_{22}[\Delta_0]=E\{\hat{\Vector{z}}_2[n]\hat{\Vector{z}}_2^T[n-\Delta_0]\}$ and $\Delta_0$ is a non-zero integer. In this new function, the correlation between the two variables $a_0[n]$ and $b_0[n]$ is not normalised by their variances, but by their own correlation for a common time lag of $\Delta_0$. $\boldsymbol{\alpha}_i$ and $\boldsymbol{\beta}_i$, $i=0,1, \dots,
q_1-1$, are all obtained in a similar way with the same normalisation. Since the noise components are not correlated with each other and not correlated with $\Vector{z}_1$ and $\Vector{z}_2$ either, we have $\hat{\boldsymbol{\Sigma}}_{12}=\boldsymbol{\Sigma}_{12}$ and for nonzero $\Delta_0$, the denominator in \eqref{eqn:gcca_firstcorrelation1} does not include any noise information. So although there is noise component existing in the original variables, the vectors $\boldsymbol{\alpha}_i$ and $\boldsymbol{\beta}_i$ obtained in this way will not depend on the noise component at all. So the effect of noise has been removed without estimating its variances.

\subsection{Applying GCCA to the BSS Problem}
\label{sec:GCCA_apply}

Applying this generalised CCA to the BSS problem, we can replace $\hat{\Vector{z}}_1[n]$ by $\Vector{x}[n]$ in \eqref{eqn:gcca_noisymixtures} and $\hat{\Vector{z}}_2[n]$ by $\Vector{x}[n-\Delta_1]$ with $\Delta_1\neq \Delta_0$. The two vectors $\boldsymbol{\alpha}_0$ and $\boldsymbol{\beta}_0$ will be the same as the extraction vector $\Vector{w}_0$~\cite{liu07l}. The extracted signal will be
\begin{equation}
y_0[n]=\Vector{w}_0^T\Vector{x}[n]\;.
\end{equation}

Then the maximization problem in \eqref{eqn:gcca_firstcorrelation} becomes
\begin{equation}
    \max_{\Vector{w}_0} \hat{J}_0(\Vector{w}_0) \;,
    \label{eqn:gcca_bssfirstcorrelation}
\end{equation}
with
\begin{equation}
\hat{J}_0(\Vector{w}_0)=\frac{E\{y_0[n]y_0[n-\Delta_1]\}}{E\{y_0[n]y_0[n-\Delta_0]\}}=\frac{\Vector{w}_0^T\Matrix{R}_{xx}[\Delta_1]\Vector{w}_0}{\Vector{w}_0^T\Matrix{R}_{xx}[\Delta_0]\Vector{w}_0}\;,
\label{eqn:gcca_bsscost}
\end{equation}
where $\Matrix{R}_{xx}[\Delta_i]=E\{\Vector{x}[n]\Vector{x}[n-\Delta_i]^T\}$, $i=0,1$,  is the correlation
matrix of the observed mixtures. From \eqref{eqn:gcca_noisymixtures}, we have
\begin{eqnarray}
\Matrix{R}_{xx}[\Delta_i]&=&\Matrix{A}E\{\Vector{s}[n]\Vector{s}^T[n-\Delta_i]\}\Matrix{A}^T\nonumber \\
&&+E\{\Vector{v}[n]\Vector{v}^T[n-\Delta_i]\}\nonumber \\
 &=&\Matrix{A}\Matrix{R}_{ss}[\Delta_i]\Matrix{A}^T\;,
\end{eqnarray}
since $\Matrix{R}_{vv}[\Delta_i]=0$ for $\Delta_i\neq 0$.

So the cost function $\hat{J}_0$ can be further simplified to
\begin{equation}
\hat{J}_0(\Vector{w}_0)=\frac{\Vector{w}_0^T\Matrix{A}\Matrix{R}_{ss}[\Delta_1]\Matrix{A}^T\Vector{w}_0}{\Vector{w}_0^T\Matrix{A}\Matrix{R}_{ss}[\Delta_0]\Matrix{A}^T\Vector{w}_0}\;.
\label{eqn:gcca_bsscost1}
\end{equation}
We assume that all of the diagonal elements of $\Matrix{R}_{ss}[\Delta_0]$ are positive, which means each of the source signals themselves should have a positive correlation value with its delayed version by $\Delta_0$.

In the next, we give a brief proof that maximization of $\hat{J}_0(\Vector{w}_0)$ with
respect to $\Vector{w}_0$ will lead to a successful extraction of
one of the source signals in the presence of noise.

Let $\Vector{g}_0=\Matrix{A}^T\Vector{w}_0$
denote the first global mixing vector. Then $\hat{J}_0(\Vector{w}_0)$ in
\eqref{eqn:gcca_bsscost1} changes into
\begin{equation}
 \hat{J}_0(\Vector{w}_0)=\frac{\Vector{g}_0^T\Matrix{R}_{ss}[\Delta_1]\Vector{g}_0}{\Vector{g}_0^T\Matrix{R}_{ss}[\Delta_0]\Vector{g}_0}\;.
\label{eqn:cca_bsscost1}
\end{equation}
Since all of the diagonal elements of the diagonal matrix $\Matrix{R}_{ss}[\Delta_0]$ are positive, we shall assume
$\Matrix{R}_{ss}[\Delta_0]=\Matrix{I}$, as the differences in the
diagonal elements of $\Matrix{R}_{ss}[\Delta_0]$ can always be absorbed
into the mixing matrix $\Matrix{A}$. This way, the diagonal
elements of $\Matrix{R}_{ss}[\Delta_1]$ become the ``normalised''
autocorrelation values of each source signal and they are assumed
to be different from each other. Note the ``normalisation'' here is not by $E\{s^2_l[n]\}$, but by $E\{s_l[n]s_l[n-\Delta_0]\}$, $l=0, 1, \dots, L-1$. Now we have
\begin{equation}
 J_0(\Vector{w}_0)=\hat{\Vector{g}}_0^T\Matrix{R}_{ss}[\Delta_l]\hat{\Vector{g}}_0\;,
 \label{eqn:cca_bsscost3}
\end{equation}
where
$\hat{\Vector{g}}_0=\frac{\Vector{g}_0}{\sqrt{\Vector{g}_0^T\Vector{g}_0}}$,
which has a property $\hat{\Vector{g}}_0^T\hat{\Vector{g}}_0=1$.

This is an eigenvalue problem and starting from here, we can use the results given in \cite{liu07l} to complete the proof and draw the conclusion that when we maximize $\hat{J}_0(\Vector{w}_0)$ with respect to $\Vector{w}_0$, this
will result in a successful extraction of the source signal with
the maximum ``normalised'' autocorrelation value.

After extracting the first source signal, we may use a deflation approach to remove it from the mixtures and then
subsequently perform the next extraction~\cite{cichocki03a}. This procedure is repeated until the last source
signal is recovered.


\section{Adaptive Realisation}
\label{sec:GCCA_adaptive}

\subsection{A direction approach}
\label{sec:GCCA_adaptive_direct}

As in the noise-free case~\cite{liu07l}, from the proof we can see that both the correlation matrix in both the numerator and the denominator in the cost function $\hat{J}_0$ can be replaced by a linear combination of the correlation matrices at different time lags to improve its robustness, as long as the one at the denominator is positive definite. More specifically, instead of maximizing the correlation between $y_0[n]$ and
$y_0[n-\Delta_1]$, we maximize the correlation between $y_0[n]$
and a weighted sum of $y_0[n-p]$, $p=2, 3, \dots, P+1$~\cite{liu08g}. Now the new
cost function is given by
\begin{equation}
J(\Vector{w}_0)=\frac{E\{y_0[n]e_0[n]\}}{E\{y_0[n]y_0[n-1]\}}\;,
\label{eqn:gcca_newcost}
\end{equation}
where
\begin{equation}
 e_0[n]=\Vector{b}^T\Vector{y}_0[n]
\end{equation}
with
\begin{eqnarray}
     \Vector{b}&=&\left[b_1 \; b_2 \; \cdots \;
     b_P\right]^\mathrm{T}\nonumber \\
\Vector{y}_0[n]&=& \left[ y_0[n-2] \; y_0[n-3] \cdots y_0[n\!-P\!-1]\right]^\mathrm{T}
\end{eqnarray}
As shown in \eqref{eqn:gcca_newcost}, we have chosen $\Delta_0=1$ because in reality more likely the signal is positively correlated with its delayed version by one sample.

Applying the standard gradient descent method to
$J_0(\Vector{w}_0,\Vector{b})$, we have
\begin{eqnarray}
\nabla_{\Vector{w}_0}J
 &=&\frac{1}{E\{y_0[n]y_0[n-1]\}^2}\left(E\{y_0[n]\hat{\Vector{x}}[n]+\right.\nonumber\\
 && \Vector{x}[n]e_0[n]\}E\{y_0[n]y_0[n-1]\}-\nonumber\\
 &&\left. E\{y_0[n]e_0[n]\}E\{y_0[n]\Vector{x}[n-1]+\Vector{x}[n]y_0[n-1]\}\right)\nonumber \\
 &&\;
\label{eqn:gcca_bssgradient}
\end{eqnarray}
where
\begin{equation}
\hat{\Vector{x}}[n]=\sum_{p=2}^{P+1}b_p\Vector{x}[n-p]\;.
\label{eqn:gcca_bssnewoutput}
\end{equation}
The correlation $E\{y_0[n]y_0[n-1]\}$ can be estimated recursively by
\begin{eqnarray}
 \sigma_{y}[n]&=&\beta\sigma_{y}[n-1]+(1-\beta)y_0[n]y_0[n-1]\;,
\label{eqn:gcca_ycorrelation}
\end{eqnarray}
where$\beta$ is the corresponding forgetting factor with $0\leq\beta<1$.

Following some standard  stochastic approximation
techniques~\cite{haykin96a},  we
obtain the following online update equation
\begin{eqnarray}
\Vector{w}_0[n+1]\!\!&=\!\!&\Vector{w}_0[n]+\nonumber\\
&&\frac{\mu}{\sigma_{y}^2}\left((y_0[n]\hat{\Vector{x}}[n]+\Vector{x}[n]e_0[n])(y_0[n]y_0[n-1])-\right.\nonumber\\
 &&\left.(y_0[n]e_0[n])(y_0[n]\Vector{x}[n-1]+\Vector{x}[n]y_0[n-1])\right)\nonumber \\
 &&\;
\label{eqn:gcca_adaptive algorithm}
\end{eqnarray}
where $\mu$ is the updating step size.

To avoid the critical case where the norm of $\Vector{w}_0[n]$
becomes too small, after each update, we normalize it to unit
length, which yields
\begin{equation}
     \Vector{w}_0[n+1] \leftarrow \Vector{w}_0[n+1]/\sqrt{\Vector{w}_0^T[n+1]\Vector{w}_0[n+1]}\;.
     \label{eqn:gcca_normaliseweight}
\end{equation}

\subsection{Adaptive Realisation Based on the Dual-Linear Predictor Structure}
\label{sec:DUAL_LP}
For noise-free mixtures, a linear predictor can be employed to extract one of the sources~\cite{liu05b,liu06c,liu07d,liu08b}, and it is closely related to the CCA approach, as shown in \cite{liu07l}. Similarly, for the GCCA approach, we can develop a corresponding dual-linear predictor structure for its implementation~\cite{liu08h}.

\subsubsection{The Structure}
\label{sec:DUAL_LPstructure}

\begin{figure}
 \begin{center}
\includegraphics[angle=0,width=0.35\textwidth]{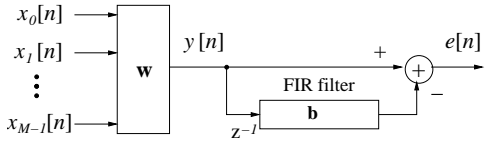}
    \caption{A linear predictor based BSE structure.
    \label{fig:dual08_LP}}
\end{center}
\end{figure}

For noise-free mixtures, a linear predictor can be employed to extract one of the sources,
as shown in Figure~\ref{fig:dual08_LP}, where the extracted signal
$y[n]$ and the instantaneous output error $e[n]$ of the linear
predictor with a length $P$ are given by
\begin{eqnarray}
y[n] &=& \Vector{w}^T\Vector{x}[n] \nonumber \\
e[n] &=& y[n]-\Vector{b}^T\Vector{y}[n]\;,
\end{eqnarray}
where $\Vector{w}$ is the demixing vector and
\begin{eqnarray}
     \Vector{b} &=& \left[b_1 \; b_2 \; \cdots \; b_P\right]^\mathrm{T} \nonumber \\
     \Vector{y}[n]&=& \left[ y[n-1] \; y[n-2]\; \cdots \; y[n\!-P]\right]^\mathrm{T}\;.
     \label{eqn:dual08_vectors}
\end{eqnarray}

The cost function is given by
\begin{equation}
 J_0(\Vector{w})=\frac{E\{e^2[n]\}}{E\{y^2[n]\}}\;.
\label{eqn:dual08_cost3}
\end{equation}
As proved in \cite{liu07d}, by minimising $J_0(\Vector{w})$
with respect to $\Vector{w}$, the sources can be extracted
successfully.

\begin{figure}
 \begin{center}
\includegraphics[angle=0,width=0.35\textwidth]{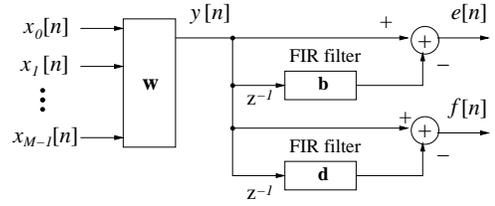}
    \caption{The proposed dual-linear predictor structure for BSE.
    \label{fig:dual08_DUALLP}}
\end{center}
\end{figure}
However, in the presence of noise, there will be a noise term in both the numerator and the denominator of \eqref{eqn:dual08_cost3} and the proof in \cite{liu07d} is not valid any more. To remove the effect of noise, as in GCCA, we propose to exploit the white nature of the noise components and employ a dual-linear predictor structure as shown in Fig.~\ref{fig:dual08_DUALLP}, where  a second linear predictor with coefficients vector $\Vector{d}$ of length $P_d$ is employed and the error signal $f[n]$ is  given by
\begin{eqnarray}
f[n] &=& y[n]-\Vector{d}^T\Vector{y}_d[n]\;,
\end{eqnarray}
where
\begin{eqnarray}
     \Vector{d} &=& \left[b_1 \; b_2 \; \cdots \; b_{P_d}\right]^\mathrm{T} \nonumber \\
     \Vector{y}_d[n]&=& \left[ y[n-1] \; y[n-2]\; \cdots \; y[n\!-P_d]\right]^\mathrm{T}\;.
     \label{eqn:dual08_dualvectors}
\end{eqnarray}

For the first linear predictor, the mean square prediction error (MSPE) $E\{e^2[n]\}$ is given by
\begin{eqnarray}
 E\{e^2[n]\}&=&E\{y^2[n]\}-2E\{y[n] \Vector{b}^T\Vector{y}[n]\}+ \nonumber \\
 && E\{\Vector{b}^T\Vector{y}[n] \Vector{y}^T[n]\Vector{b}\} \nonumber \\
   &=&\sum_{p=0}^{P}b_p^2E\{y^2[n-p]\}-\nonumber\\
   &&\sum_{p,q=0;p\neq q}^{P}s_{pq}b_pb_qE\{y[n-p]y[n-q]\}\nonumber \\
   &=&q_cE\{y^2[n]\}-\nonumber \\
   &&\sum_{p,q=0;p\neq q}^{P}s_{pq}
   b_pb_q\Vector{w}^T\Matrix{R}_{xx}[q-p]\Vector{w}\nonumber \\
   &=&q_cE\{y^2[n]\}-\nonumber\\
   &&\Vector{w}^T(\sum_{p,q=0;p\neq q}^{P}s_{pq}
   b_pb_q\Matrix{R}_{xx}[q-p])\Vector{w}\;,
\label{eqn:dual08_mspee}
\end{eqnarray}
where  $q_c=\sum_{p=0}^{P}b_p^2$ with $b_0=1$, and $s_{pq}$ is $1$ when $p=0$ or $q=0$, and $-1$ otherwise.
From \eqref{eqn:gcca_noisymixtures}, \eqref{eqn:gcca_noisecorrelation} and \eqref{eqn:gcca_sourcecorrelation}, we have
\begin{eqnarray}
\Matrix{R}_{xx}[p-q]&=&\Matrix{A}E\{\Vector{s}[n]\Vector{s}^T[n-(p-q)]\}\Matrix{A}^T\nonumber \\
&&+E\{\Vector{v}[n]\Vector{v}^T[n-(p-q)]\}\nonumber \\
 &=&\Matrix{A}\Matrix{R}_{ss}[p-q]\Matrix{A}^T\;,
\end{eqnarray}
for $p\neq q$. Then we have
\begin{eqnarray}
 E\{e^2[n]\}&=&q_cE\{y^2[n]\}-\nonumber\\
 &&\Vector{w}^T\Matrix{A}(\sum_{p,q=0;p\neq q}^{P}s_{pq}
   b_pb_q\Matrix{R}_{ss}[q-p])\Matrix{A}^T\Vector{w}\nonumber \\
   &=&q_cE\{y^2[n]\}-\nonumber\\
   &&\Vector{g}^T(\sum_{p,q=0;p\neq q}^{P}s_{pq}
   b_pb_q\Matrix{R}_{ss}[q-p])\Vector{g}\nonumber\\
   &=&q_cE\{y^2[n]\}-\Vector{g}^T\hat{\Matrix{R}}_{ss}\Vector{g}\;,
\label{eqn:dual08_mspee1}
\end{eqnarray}
with $\Vector{g}=\Matrix{A}^T\Vector{w}$ denoting the global demixing
vector and  $\hat{\Matrix{R}}_{ss}$ is a diagonal matrix given by
\begin{equation}
\hat{\Matrix{R}}_{ss}=\sum_{p,q=0;p\neq q}^{P}s_{pq}
   b_pb_q\Matrix{R}_{ss}[q-p]\;,
\end{equation}
with its $l-th$ diagonal element $\hat{r}_l$ given by
\begin{equation}
\hat{r}_l=\sum_{p,q=0;p\neq q}^{P}s_{pq}
   b_pb_q\rho_l[q-p]\;.
\end{equation}

Similarly, for the second linear predictor, we have
\begin{eqnarray}
 E\{f^2[n]\}&=&E\{y^2[n]\}-2E\{y[n] \Vector{d}^T\Vector{y}_d[n]\}+ \nonumber \\
 && E\{\Vector{d}^T\Vector{y}_d[n] \Vector{y}_d^T[n]\Vector{d}\} \nonumber \\
   &=&a_cE\{y^2[n]\}-\nonumber \\
   &&\Vector{g}^T(\sum_{p,q=0;p\neq q}^{P_d}s_{pq}
   d_pd_q\Matrix{R}_{ss}[q-p])\Vector{g}\nonumber \\
   &=&a_cE\{y^2[n]\}-\Vector{g}^T\tilde{\Matrix{R}}_{ss}\Vector{g}\;,
\label{eqn:dual08_mspef}
\end{eqnarray}
with $a_c=\sum_{p=0}^{P_d}d_p^2$ with $d_0=1$ and
 $\tilde{\Matrix{R}}_{ss}$ is a diagonal matrix given by
\begin{equation}
\tilde{\Matrix{R}}_{ss}=\sum_{p,q=0;p\neq q}^{P_d}s_{pq}
   d_pd_q\Matrix{R}_{ss}[q-p]\;,
\end{equation}
with its $l-th$ diagonal element $\tilde{r}_l$ given by
\begin{equation}
\tilde{r}_l=\sum_{p,q=0;p\neq q}^{P_d}s_{pq}
   d_pd_q\rho_l[q-p]\;.
\end{equation}

\subsubsection{The Proposed Cost Function}
\label{sec:DUAL_LPcost}

Note in the second term of both \eqref{eqn:dual08_mspee1} and \eqref{eqn:dual08_mspef}, there is not any noise component. Then we can construct a new cost function as follows
\begin{equation}
J(\Vector{w})=\frac{q_cE\{y^2[n]\}-E\{e^2[n]\}}{a_cE\{y^2[n]\}-E\{f^2[n]\}}=\frac{\Vector{g}^T\hat{\Matrix{R}}_{ss}\Vector{g}}{\Vector{g}^T\tilde{\Matrix{R}}_{ss}\Vector{g}}\;.
\label{eqn:dual08_cost}
\end{equation}
Now we impose another condition on the second linear predictor: suppose the coefficients $\Vector{d}$ are chosen in such a way that all of the diagonal elements of $\tilde{\Matrix{R}}_{ss}$ are of positive value. This is a difficult condition due to the blind nature of the problem. However, for a special case with $P_d=1$ and $d_1=1$, i.e. a one step ahead predictor, we have
\begin{equation}
\tilde{\Matrix{R}}_{ss}=2\Matrix{R}_{ss}[1]\;,
\end{equation}
which is the correlation matrix of the source signals with a time lag of $1$. Then the condition means each of the source signals should have a positive correlation with a delayed version of itself by lag $1$. As mentioned in Sec.~\ref{sec:GCCA_adaptive_direct}, in reality, there are many signals having this correlation property and therefore can meet this requirement. Now we can see the cost function has the same form as in \ref{eqn:cca_bsscost1}. Therefore, we can consider this dual-linear predictor structure as an indirect implementation of the GCCA approach for solving the BSS problem.

Since all of the diagonal elements of $\tilde{\Matrix{R}}_{ss}$ are positive, we shall assume
$\tilde{\Matrix{R}}_{ss}=\Matrix{I}$, i.e. $\tilde{r}_l=1$, $l=0, 1, \dots, L-1$, as the differences in the
diagonal elements can always be absorbed into the mixing matrix $\Matrix{A}$. This way, the diagonal
elements $\hat{r}_l$, $l=0, 1, \dots, L-1$, of $\hat{\Matrix{R}}_{ss}$ in the numerator become the ``normalised''
autocorrelation values of each source signal and they are assumed
to be different from each other. For the case with $P_d=1$ and $d_1=1$, the ``normalisation'' here is not by $E\{s^2_l[n]\}$, but by $\tilde{r}_l=2E\{s_l[n]s_l[n-1]\}$.

 Now we have
\begin{equation}
 J(\Vector{w})=\hat{\Vector{g}}^T\hat{\Matrix{R}}_{ss}\hat{\Vector{g}}\;,
 \label{eqn:dual08_cost2}
\end{equation}
where
$\hat{\Vector{g}}=\frac{\Vector{g}}{\sqrt{\Vector{g}^T\Vector{g}}}$,
which has a property $\hat{\Vector{g}}^T\hat{\Vector{g}}=1$. Clearly, according to the proof provided earlier, we can draw the conclusion that when we minimize $J(\Vector{w})$ with respect to $\Vector{w}$, this
will result in successful extraction of the source signal with
the minimum ``normalised'' autocorrelation value. Note here the extracted signal is not the one with the maximum ``normalised'' autocorrelation value.


\subsubsection{Adaptive Algorithm}
\label{sec:DUAL_LPalgorithm}

Applying the standard gradient descent method to
$J(\Vector{w})$, we have
\begin{eqnarray}
 \nabla_{\Vector{w}}J&=&\frac{2}{(a_cE\{y^2[n]\}-E\{f^2[n]\})^2}\left((q_cE\{y[n]\Vector{x}[n]\}-\right.\nonumber\\
 &&  E\{e[n]\hat{\Vector{x}}[n]\})(a_cE\{y^2[n]\}-E\{f^2[n]\})-\nonumber\\
 && (q_cE\{y^2[n]\}-E\{e^2[n]\})(a_cE\{y[n]\Vector{x}[n]\}-\nonumber\\
 &&\left.E\{f[n]\tilde{\Vector{x}}[n]\})\right)\;,
\label{eqn:dual08_gradient}
\end{eqnarray}
where
\begin{eqnarray}
\hat{\Vector{x}}[n]&=&\Vector{x}[n]-\sum_{p=1}^{P}b_p\Vector{x}[n-p]\nonumber \\
\tilde{\Vector{x}}[n]&=&\Vector{x}[n]-\sum_{p=1}^{P_d}d_p\Vector{x}[n-p]\;.
\label{eqn:dual08_newoutput}
\end{eqnarray}
$E\{e^2[n]\}$, $E\{y^2[n]\}$ and $E\{f^2[n]\}$ can be estimated respectively by
\begin{eqnarray}
 \sigma_e[n]&=&\beta_e\sigma_e[n-1]+(1-\beta_e)e^2[n]\nonumber \;,\\
 \sigma_y[n]&=&\beta_y\sigma_y[n-1]+(1-\beta_y)y^2[n]\nonumber \;,\\
 \sigma_f[n]&=&\beta_f\sigma_f[n-1]+(1-\beta_u)f^2[n]\;,
\label{eqn:dual08_powerest}
\end{eqnarray}
where $\beta_e$, $\beta_y$ and $\beta_{f}$ are the corresponding
forgetting factors with $0\leq\beta_e, \beta_y, \beta_{f}<1$.

Following standard  stochastic approximation
techniques~\cite{haykin96a}, we
obtain the following online update for $\Vector{w}[n]$
\begin{eqnarray}
\Vector{w}[n+1]&=&\Vector{w}[n]-\frac{2\mu}{(a_c\sigma_y-\sigma_f)^2}\left((q_cy[n]\Vector{x}[n]-\right.\nonumber \\
&&\left. e[n]\hat{\Vector{x}}[n])(a_c\sigma_y-\sigma_f)\right.-(q_c\sigma_y-\sigma_e)\cdot\nonumber\\
&&\left.(a_cy[n]\Vector{x}[n]-f[n]\tilde{\Vector{x}}[n])\right)\;,
\label{eqn:dual08_updatew}
\end{eqnarray}
where $\mu$ is the learning rate.
For the case with $P_d=1$ and $d_1=1$, we have $a_c=2$ in \eqref{eqn:dual08_updatew}, which will be used in our simulations.

%
%
\section{Simulations}
\label{sec:GCCA_sim}

\begin{figure}
\begin{center}
   \includegraphics[angle=0,width=0.45\textwidth]{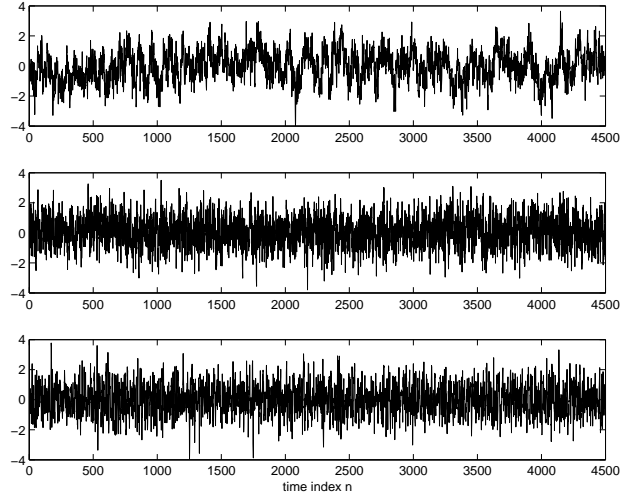}
   \caption{The three source signals used in the simulations.
    \label{fig:dual08_sources}}
\end{center}
\end{figure}
Here we only provide some preliminary simulation results based on the dual-linear predictor structure~\cite{liu08h}. Three source signals are used which are generated by passing three randomly generated white Gaussian signals through three different filters. The power of the sources is normalised to one. The correlation value of each of the source signals is checked to make sure it is positive and not close to zero for one sample shift. Fig.~\ref{fig:dual08_sources} shows the three source signals,
denoted by $s_0$, $s_1$ and  $s_2$, respectively.

The coefficients of the first linear predictor coefficients $\Vector{b}$  were randomly generated with a length of $P=5$, and given by
\begin{equation}
     \Vector{b}=[-0.4548\;   -1.0053\;    1.1957\;   -0.5590\;   -0.3617]\;.
     \label{eqn:dual08_predictorcoefficients}
\end{equation}
For the second linear predictor, $P_d=1$, $d_1=1$, and $a_c=2$.

The normalised correlation value $\hat{r}_l$ for each source signal with this dual-linear predictor configuration is $0.0395$,   $0.2174$ and   $0.7949$, respectively. As already proved, since the first source signal has the smallest correlation value of $0.0395$, it will be extracted by minimizing the cost function.

\begin{figure}
\begin{center}
   \includegraphics[angle=0,width=0.45\textwidth]{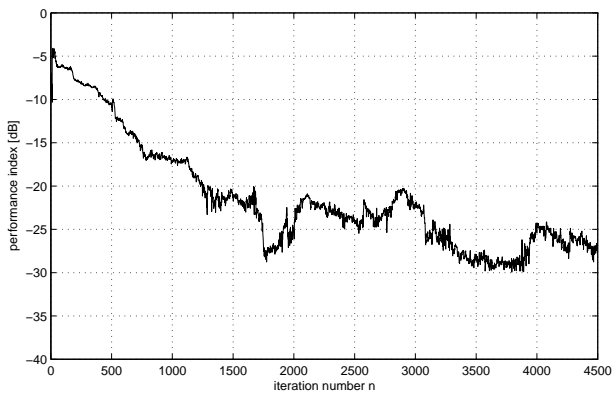}
   \caption{One of the learning curves obtained in our simulations.
    \label{fig:dual08_onecurve}}
\end{center}
\end{figure}

The $3\times3$ mixing matrix $\Matrix{A}$ is randomly generated
and given by
\begin{equation}
  \Matrix{A}= \left[ \begin{array}{ccc}
  0.9207   & 0.0299  &  0.3891\\
    0.5165 &   0.3676  &  0.7733\\
    0.7822 &  -0.2735 &  -0.5598  \end{array} \right]\;.
\end{equation}
Its row vector is normalised to unity to make sure it is comparable to the noise variance, which is  $\sigma_v^2=0.09$.  The forgetting factors is $\beta_e=\beta_y=\beta_f=0.975$ and the stepsize
$\mu=0.0015$. A learning curve for this case is shown in
Fig.~\ref{fig:dual08_onecurve}, with the performance index
defined as~\cite{cichocki03a}
\begin{equation}
PI=10\log_{10}\left(\frac{1}{L-1}(\sum_{l=0}^{L-1}\frac{g_l^2}{\max\{g_0^2, \dots, g_{L-1}^2\}}-1)\right), \label{eqn:pi_bse2}
\end{equation}
with $\Vector{g}=[g_0\; g_1 \;\cdots\; g_{L-1}]$.

\begin{figure}[t]
\begin{center}
   \includegraphics[angle=0,width=0.45\textwidth]{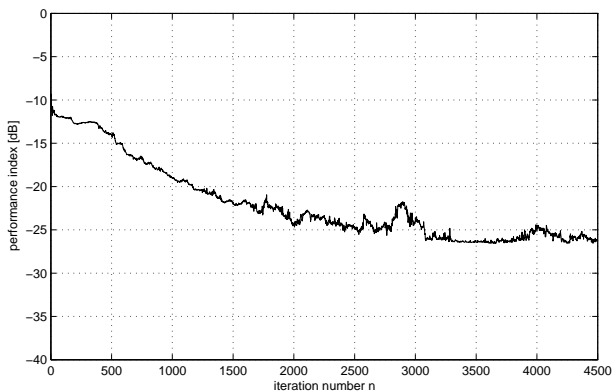}
   \caption{The averaged learning curve over $1000$ runs.
    \label{fig:dual08_learningcurve}}
\end{center}
\end{figure}

To show its performance in a more general context, we change the initial value of the demixing vector $\Vector{w}$ randomly each time to run the algorithm and the average learning curve over $1000$ runs is given in Fig.~\ref{fig:dual08_learningcurve}.  Both curves show a successful extraction of the source signal.

%
%
\section{CONCLUSIONS}
\label{sec:GCCA_concl} The traditional canonical correlation analysis has been generalised to include noisy signals where the effect of noise can be eliminated effectively by the proposed approach. It was then applied to the blind source separation problem and adaptive implementations were derived. In particular, a dual-linear predictor structure was proposed to blindly extract the source signals from their noisy mixtures, and it can be considered as an indirect implementation of GCCA. Some preliminary simulation results have been provided to show the effectiveness of the proposed approach.


\end{document}